\documentclass{ifacconf}

\usepackage{natbib,graphics, graphicx, amsmath, float, amsfonts, amssymb, fp, mathrsfs, booktabs} 

\counterwithin*{section}{part}

\usepackage{thmtools}
\usepackage{cleveref}
\makeatletter
\newcommand{\labelx}[1]{
    \relax
    \ifmmode
        \label{#1} 
    \else 
        \ifnum\pdfstrcmp{\@currenvir}{document}=0
            \label{#1}
        \else
            \label[\@currenvir]{#1}
        \fi
    \fi
}
\makeatother

\newtheorem{proofpart}{Part}
\newtheorem{proofcase}{Case}

\renewcommand{\vec}[1]{\mathbf{#1}}
\newcommand{\bi}{\begin{itemize}}
\newcommand{\ei}{\end{itemize}}
\newcommand{\beq}{\begin{enumerate}}
\newcommand{\eeq}{\end{enumerate}}
\renewcommand{\vec}[1]{\mathbf{#1}}
\newcommand{\be}{\begin{equation}}
\newcommand{\ee}{\end{equation}}

\newcommand{\bc}{\begin{center}}
\newcommand{\ec}{\end{center}}
\newcommand{\bd}{\begin{defn}}
\newcommand{\ed}{\end{defn}}
\newcommand{\bt}{\begin{thm}}
\newcommand{\et}{\end{thm}}
\newcommand{\bp}{\begin{proof}}
\newcommand{\ep}{\end{proof}}
\renewcommand{\l}{\left}
\renewcommand{\r}{\right}

\newcommand{\R}{\mathbb{R}}
\newcommand{\V}{\vert}

\DeclareMathOperator*{\argmin}{arg\,min}

\crefname{thm}{Theorem}{Theorems}
\crefname{cor}{corollary}{corollaries}
\crefname{lem}{lemma}{lemmas}
\crefname{prop}{proposition}{propositions}
\crefname{con}{conjecture}{conjectures}
\crefname{defn}{definition}{definitions}
\crefname{eg}{example}{examples}
\crefname{rem}{remark}{remarks}
\crefname{nota}{notation}{Notation}
\crefname{assum}{assumption}{Assumption}
\crefname{ceg}{Counter-example}{Counter-examples}
\crefname{equation}{}{}
\begin{document}
\begin{frontmatter}

\title{Optimal Responses to Constrained Bolus Inputs to Models of T1D}


\author[First]{Christopher Townsend}
\author[First]{Maria M. Seron}
\author[Second]{Nicolas Magdelaine}

\address[First]{School of Engineering, University of Newcastle, Australia (emails: chris.townsend@newcastle.edu.au, maria.seron@newcastle.edu.au)}

\address[Second]{L@bISEN, robotics, ISEN Yncr\'ea Ouest, France (email: nicolas.magdelaine@isen-ouest.yncrea.fr)}

\begin{abstract}
We characterise the bolus insulin  input which minimises the maximum plasma glucose concentration predicted by the Magdelaine and Bergman minimal models in response to any positive bounded disturbance whilst remaining above a fixed lower plasma glucose concentration. This characterisation is in terms of the maxima and minima of the plasma glucose concentration and limits the controllability of such systems. Any further attempt to lower the maximum plasma glucose concentration will result in hypoglycaemia. 
\end{abstract}

\begin{keyword}
Optimal control, Model-based Predictive Control
\end{keyword}

\end{frontmatter}

\section{Introduction}

Type one diabetics are unable to regulate plasma glucose levels which if not successfully controlled result in several adverse health outcomes. Diabetes is a chronic, life-long disease affecting over thirty-eight million people \citep{wyou16}. Currently, a diabetic's plasma glucose concentration is controlled by the subcutaneuous administration of insulin to minimise plasma glucose concentrations whilst avoiding hypoglycaemia. Insulin requirements vary depending on a variety of physiological factors and external disturbances. Thus to improve control and reduce the burden of management, research efforts have been focussed on the development of an artificial pancreas \citep{harv10}.

Models of the glucose insulin dynamics in type one diabetics assist in the development of such systems and current management for example by predicting future glucose concentrations based on current inputs. A number of models of glucose regulation have been proposed (\cite{makr06, wili09, colm14}). Each is typically comprised of sub-systems describing different physiological processes such as insulin kinetics and glucose absorption.

 Recently, research has focused on comprehensive models of glucose dynamics. Typically, these models are high order non-linear dynamic systems with many parameters to ensure robustness to inter-individual variability. However, simpler models are useful to establish general theoretical properties that would otherwise be difficult to investigate analytically. Indeed, most models of glucose dynamics share certain analytic properties -- such as positivity of the plasma glucose. Thus analytic results obtained for simpler models can give insights into the behaviour of more comprehensive models. Hence, we focus on analytic properties of the \emph{Magdelaine} \citep{gonz17,magd15,riva17} and \emph{Bergman models} \citep{kand09,berg05} of glucose-insulin dynamics.

The need to avoid the hypoglycaemic threshold whilst minimising the maximum glucose concentration whilst the system is subject to bounded disturbances means that the control of blood glucose concentrations may be considered as a constrained optimisation problem.

The work of \cite{town17,town17IFAC} presented fundamental control limitation for the minimisation of the maximum glucose concentration in the Bergman Minimal model \cite{berg05} when the bolus insulin input was constrained to be a pulse input. It was proven that if the maxima and minima of the glucose concentration are interlaced then the maximum glucose concentration is minimised and any attempt to further lower this maximum will result in hypoglycaemia.

In this paper we extend this characterisation to the Magdelaine model which imposes an additional constraint on the insulin input $u$.

Insulin inputs are usually separated into \emph{basal} inputs which are, typically constant inputs, used to keep the system in equilibrium in the absence of exogenous disturbances and \emph{bolus} inputs which are bounded inputs delivered to move the system from equilibrium or minimise the impact of exogenous disturbances. So for a model of plasma glucose concentration with insulin input $u(t)$ and output $g$ which represents the plasma glucose concentration, the insulin input $u$ is a positive real function of the form:
\begin{align}
\labelx{eq:u}
	u(t) =\hat u(t) + \overline u(t)
	\end{align}%
where $\hat u(t)$ is the \emph{bolus} input and $\overline u(t)$ is the \emph{basal} input. Additionally, the basal input, $\overline u (t)$, is such that the first derivative of the response, $g$, satisfies $\dot g = 0$ in the absence of exogenous disturbances and, if possible, the steady-state glucose concentration, $g(\infty)$, equals a specified concentration. As will be explored here, in the Magdelaine model the basal input can only achieve the first criterion, as the derivative of the plasma glucose concentration is independent of the current concentration meaning that the steady-state glucose concentration is independent of the basal input. In contrast, in the Bergman model the steady-state glucose concentration is a globally asymptotically stable equilibrium determined by the basal input $\overline u$. Thus both criteria may be met simultaneously. 

The total bolus insulin was not constrained by \cite{town17,town17IFAC}. However, as we require the plasma glucose concentration to return to a specified concentration and the steady-state glucose concentration in the Magdelaine model is an unstable equilibrium, we characterise the optimal input when the total bolus input is constrained. We then consider optimality of such constrained inputs to the Bergman model.

Here we do not propose a control strategy but rigourously prove the limitations in the controllability of the Magdelaine and Bergman models subject to bounded disturbances. An exploration of the clinical implications of such control limitations is given in \cite{town18}. Furthermore, we believe, a mathematical and rigourous understanding of the models of type one diabetes is necessary for the development of controllers based on such models.



\section{Magdelaine Model and Constraints}
The \emph{Magdelaine Model} is the affine system:
\begin{align}
\labelx{eq:magdmat}
\begin{split}
	\dot x_1  & = -\alpha_2 x_2 + \alpha_4 x_4 + E\\
	\dot x_2  & = -\alpha_3 x_2 + \alpha_3 x_3\\
	\dot x_3  & = -\alpha_3 x_3 +\alpha_3 u\\
	\dot x_4  & = -\alpha_5 x_4 + \alpha_5 x_5\\
	\dot x_5  & = -\alpha_5 x_5 + \alpha_5 d
	\end{split}
\end{align}%
where $u$ is the insulin input, $d$ is some positive bounded disturbance, $E$ is the endogenous glucose production and $\alpha_i \in \R _+$ are constants. The co-ordinates $x_1, x_2$ and $x_4$, of the state $x$, represent the plasma glucose, insulin effectiveness and the impact of the disturbance $d$. The states $x_3$ and $x_5$ are the subcutaneous and absorption transitional compartments. For notational simplicity we denote by $g:=x_1$, $x:=x_2$ and $w:=x_4 + E$. Also $a := \alpha _2$ and $b := \alpha_4$.

Aside from the positivity and boundedness assumptions, the disturbance $d$ is unconstrained. As outlined in \eqref{eq:upulse} the input $u$ is constrained to be a single pulse input of some finite duration.
 
 We normalise $E$ by the constant $b$ that is, $E$ in the first equation of \eqref{eq:magdmat} is replaced by $b^{-1} E$. Thus by \eqref{eq:magdmat} the plasma glucose $g$ is the solution to the differential equation:
\begin{align}
\labelx{eq:madgnew}
	\dot g = -  ax + bw
\end{align}%
	where $x$ is the insulin effectiveness and $w$ combines the endogenous glucose production and the response, $x_4$, to a positive disturbance $d$. In the absence of any disturbance we see that:
	\[
		\dot g = -ax + bE
		\]%
	We assume the bolus input has compact support. Thus the plasma glucose is in steady-state i.e. $\dot g = 0$ if and only if:
	\[
		x = \l ( \frac{b}{a} \r ) E
		\]%
	Thus for the response $g$ to be bounded it requires the input $\overline u = \frac{b}{a} E$. Therefore, as $\overline u$ is uniquely determined by $E$, we may consider the equivalent system:
\begin{align}
	\labelx{eq:model}
	\dot g = -a x + b w 
\end{align}%
where the basal input $\overline u = 0$, $E = 0$ and the set point $g(0) = 0$. We note that this is not physiological. However, setting $E = 0$ does not affect the dynamics of the system as with $\overline u$ determined as above the physiological system is a scalar offset of the system \eqref{eq:model}.

After any disturbance we require that the system return to its set point i.e. $\lim_{t \to \infty} g(t) = 0$. As the solution, to \eqref{eq:model} is:
\begin{align}
\label{eq:anumberforyou}
	g(t) = - a\int_0 ^t x \, dt  + b\int_0 ^t w \, dt
	\end{align}%
the magnitude of $x$ is bounded by the magnitude of $w$. Indeed, by \eqref{eq:magdmat} and \eqref{eq:anumberforyou}:
\begin{align}
\labelx{eq:whyareyoustillnumberingthese}
	\int_{\R _+} u(t) \, dt = \int_{\R_+} x(t) \, dt = \l ( \frac{b}{a} \r ) \int_{\R _+} w(t) \, dt 
	\end{align}%
As the system is required to return to steady-state, we have that $g(\infty) = 0$,  $x_2 (\infty) = x_2 (0)$ and $x_3 (\infty) = x_3(0)$. Thus the first equality of \eqref{eq:whyareyoustillnumberingthese} is established by integrating the second and third equations of \eqref{eq:magdmat} and the second by rearranging \eqref{eq:anumberforyou}. The equality between the volume of the bolus $u$ and the disturbance $w$ given by \eqref{eq:whyareyoustillnumberingthese} motivates \Cref{defn:adua}.

\begin{defn}[Adequate]
\labelx{defn:adua}
	Let $u(t)$ be an input and let $U$ be the \emph{amount} (1-norm) of $u(t)$:
	\[
		U:= \int_{\mathbb{R_+}} u(t) \, dt
	\]%
	$U$, or $u(t)$, is \emph{adequate} if $\lim_{t \to \infty} g(t) = 0$.
\end{defn}

We assume throughout that all inputs to the Magdelaine model are adequate.

As in \cite{town17} and \cite{town17IFAC}, we desire that there exists a fixed lower bound $\lambda$ such that $g(t) \geq \lambda$ for all $t$. We also require the function $d(t)$ to be positive, bounded and such that there exists a solution to \eqref{eq:model} and that: 
\[
	\int_{\R _+} d(t)
\]%
is bounded. We will see that the optimality conditions for the Magdelaine model are similar but not identical to those derived in \cite{town17} and \cite{town17IFAC}. This is as we require the system to return to steady-state. Without this constraint the results of \cite{town17} and \cite{town17IFAC} apply directly. Furthermore, if the optimality conditions of \cite{town17} and \cite{town17IFAC} are met by the response of the Magdelaine model rather than the conditions proposed here, then the maximum of the response will in general be lower.

\labelx{sec:5.2}

\section{Response to Pulse Inputs}


As in \cite{town17} and \cite{town17IFAC} we consider the response of the Magdelaine model to pulse inputs $u$, \eqref{eq:u}, of the form:
\begin{align}
	\labelx{eq:upulse}%
	u(t) = \overline u + \hat u \chi_A
\end{align}%
where $\overline u, \hat u \in \mathbb{R}_+$ are the \emph{basal} input and magnitude of the \emph{bolus} input respectively, and $\chi_A$ is the indicator function over a compact interval $A$. As mentioned above we may assume $\overline u  = 0$.



We constrain the response $g(t)$ by requiring that there is a fixed lower bound, $\lambda$, at or above the hypoglycaemic threshold, such that $g(t) \geq \lambda$ for all $t$.

\begin{defn}[$\lambda$--incident]
	An input $u$ is \emph{$\lambda$--incident} if the response $g(u) \geq \lambda$ for all $t$ and there exists $t_{\min}$ at which $g(u, t_{\min}) = \lambda$.
\end{defn}

	When $\lambda$ is unambiguous, we say $u$ is \emph{incident}. As we assume inputs $u$ are adequate an input can only be incident if:
	\begin{align}
	\label{eq:lambdabound}
		\lambda \geq - b \int_0 ^\infty w
	\end{align}%
	This is as the system must return to steady-state which bounds the magnitude of $u$ by the disturbance $w$, shown by \eqref{eq:whyareyoustillnumberingthese}. We fix $w$ and choose $\lambda$ such that there is an adequate, $\lambda$-incident $u$ i.e. the lower bound $\lambda$ is achievable.

	\Cref{lem:adequate} proves the existence of incident adequate inputs for any fixed lower bound and bounded disturbance $w$. In \Cref{lem:adequate} the input time of the bolus $\hat u \chi_A$ is denoted by $t'$ and the duration by $\tau$ i.e. the input $u(t) = \hat u$ for all $t \in [t', t' + \tau] =: A$.

\begin{lem}[Adequate and Incident Input] 
	\labelx{lem:adequate}
	For any $w$ and $\lambda \leq g(0)$, there exists an adequate input $u$ of the form $\hat u \chi_A$. Furthermore, if we let the input time be any real number $t' \in \R$, then there exist $t'$ and $\tau$ such that $u$ is incident.
\end{lem}

\begin{pf}
	Fix $w$ and $\lambda \geq -b M$. A solution for $g$ is:
	\[
		g(t) = -a \int_0 ^t x + b \int_0 ^t w
	\]%
	We have assumed that the norm of the disturbance $d(t)$ is bounded i.e. there exists $M \in \mathbb{R}_+$ such that:
	\[
		\int_0 ^\infty d(t) = M
	\]%
	This implies, by the fourth equation in \eqref{eq:magdmat}:
	\[
		\int_0 ^\infty w(t) = M
	\]%
	Independently of the input time and duration, $t'$ and $\tau$, there exists $\hat u$ such that:
	\[
		\int_0 ^\infty x(t) = \l (\frac{b}{a} \r )M
	\]%
	Thus:
	\[
		\lim_{t \to \infty} g(t) = 0
	\]%
	i.e. $u$ is adequate. Fix $\tau \geq 0$ and take $T > 0$. The point $T$ is arbitrary and chosen to provide sufficient time for $g$ to decrease before a positive disturbance occurs. The value of $T$ represents a prebolus interval and depends on the constants $\alpha_i$ in \eqref{eq:magdmat}.

	Suppose $d(t) = 0$ for all $t \leq T$. Then for $g(T) = \lambda$ we require:
	\[
		\int_ {0} ^T x (t) = -\frac{\lambda}{a}
	\]%
	for sufficiently large $T$ there will always exist such $u$. For $g(t) \geq \lambda$ for all $t$ we require:
	\[
		\int_{0} ^t x(t)  \leq \l ( \frac{1}{a} \r ) \l ( b \int_{0} ^t w - \lambda \r ) 
	\]%
	for all $t$. For each $\varepsilon \in (0, |\lambda|)$ there exists $T$ such that:
	\[
		\int_0 ^T w = M - \varepsilon
	\]%
	Choosing $t' > T$ ensures that:
	\[
		g(T) = b(M - \varepsilon) \geq b(M + \lambda)
	\]%
	there is no $x$--component as $u (t)  = 0$ for all $t< T< t'$. Thus applying adequate $u$ with input time $t' > T$ we have that $g(t) \geq \lambda$ for all $t$. Finally as $g$ is a continuous function of $u$, $t'$ and $\tau$ there exists an incident input $u$.
\end{pf}

The comparison of the response to distinct inputs $u$ and $v$ to characterise the response with the lowest maximum relies on the location of the intersection points of the responses $g(u)$ and $g(v)$. We later prove that the maximum of the response $g(u)$ is monotonic when the sequence of inputs are nested, see \Cref{defn:nest}. 

\begin{defn}[Nested]
\label{defn:nest}
	Suppose $u$ and $v$ are two pulses with input times $t'$ and $s'$ and durations $\tau$ and $\sigma$ respectively. Then $u$ is \emph{nested} in $v$ if $[t', t'+\tau] \subset [s', s'+ \sigma]$.
\end{defn}

\begin{lem}[Intersection Points]
	\labelx{lem:intpoints}
	Suppose $u$ and $v$ are distinct inputs to \eqref{eq:magdmat}. Then for all solutions, $\phi$, there exist at most two $t_i$ such that $\phi(t_i, u) = \phi(t_i, v)$ and these $t_i$ are distinct if and only if $u$ and $v$ are nested.
\end{lem}

\begin{pf}
	Let the input times of $u$ and $v$ be $t'$ and $s'$ respectively. Observe that $y(u) = y(v)$ if and only if $u - v$ changes sign. As $u$ and $v$ are rectangular $u-v$ can change sign at most twice. We proceed similarly for $x$.
\end{pf}

\Cref{lem:integrals} applies \Cref{lem:intpoints} to the case of the Magdelaine model.

\begin{lem}
	\labelx{lem:integrals}
	Suppose $u$ and $v$ are continuous functions which intersect $n$ times. Then for each solution to the differential equations:
	\[
		\dot x_i = - a x_i + a x_{i-1}
	\]%
	where $a>0, x_i = x_i(u)$ and $x_0 (u) = u$, the functions:
	\[
		\int_0 ^t x_i(u), \quad \text{and } \int_0 ^t x_i(v)
	\]%
	intersect at most $n-1$ times.
\end{lem}

\begin{pf}
	By \Cref{lem:intpoints} we see that if $u$ and $v$ intersect $n$ times. Then $x_i(u)$ and $x_i(v)$ may intersect at most $n$ times.  We also observe that should $t_{i,k}$ be an intersection point of $x_i(u)$ and $x_i(v)$. Then  $s_{i,k}$, the intersection point of:
	\[
		\int_0 ^t x_i(u), \quad \text{and } \int_0 ^t x_i(v)
	\]%
	resulting from $t_{i,k}$ must satisfy $s_{i,k} > t_{i,k}$. We now proceed by induction. 
	Let $t_n$ denote the $n^{\text{th}}$ intersection point of $x_1(u)$ and $x_1(v)$ after which we assume without loss of generality that $x_1(u) > x_1(v)$. Suppose there exists $s_n > t_n$ which is the $n^{\text{th}}$ intersection point of:
	\[
		\int_0 ^t x_1(u), \quad \text{and } \int_0 ^t x_1(v)
	\]%
\quad \pagebreak

Thus:
\begin{align*}
	x_1(s_n, u) & = - a \int_0 ^{s_n} x_1(u) + a \int_0 ^{s_n} u \\
	            & \geq - a \int_0 ^{s_n} x_1(v) + a \int_0 ^{s_n} v \\
		    & = x_1(s_n, v)
\end{align*}%
This implies that $s_n > t_n > t_k$, for all $k < n$ is an additional intersection point of $x_1(u)$ and $x_1(v)$ contradicting the fact that they intersect at most $n$ times. The proof now follows by induction.
\end{pf}

\section{Optimal Inputs}

We say an input is optimal if it minimises the maximum of the response compared to all other inputs whilst meeting the constraints. This is formalised in \Cref{defn:mini}. We notate the maximum of a response to an input $u$ by $\gamma(u)$ i.e. given an input $u$ we define $\gamma(u) := \max \{ g(t) \}$.

\begin{defn}[Optimal]
\label{defn:mini}
	For fixed $w$ and $\lambda$ a response $g$ is \emph{minimised} by an input $u$ if $\gamma(u) \leq \gamma(v)$ for all $v \neq u$. In which case $u$ is \emph{optimal}.
\end{defn}

\Cref{lem:youshouldhavenameditbynow} proves that the lower the fixed minimum the lower the maximum glucose concentration.

\begin{lem}
\labelx{lem:youshouldhavenameditbynow}
	Suppose either $t'$ or $\tau$ is fixed. Then the maximum $\gamma$ is a monotone function of the minimum $\lambda$.
\end{lem}

\begin{pf}
	Take $\lambda' > \lambda$ and suppose $u$ is an input which is $\lambda$--incident and $v$ is a $\lambda'$--incident input. As either $t'$ or $\tau$ are fixed $u$ and $v$ are not nested. Hence they intersect at most once. Thus by \Cref{lem:intpoints,lem:integrals} there exists no $t_g > \min\{t', s'\}$ where $t'$ and $s'$ are the respective delivery times for $u$ and $v$, such that $g(u) = g(v)$.
	As $\lambda < \lambda'$ and $u$ and $v$ are $\lambda$--incident and $\lambda'$--incident respectively, $g(u) < g(v)$ for all $t > \min\{t',s'\}$. Otherwise there would exist $t_g$ such that $g(u) = g(v)$. Thus $\gamma(u) < \gamma(v)$.
\end{pf}

An interesting property of the Magdelaine model is that the minimisation of the maximum of the response $g(t)$ is equivalent to minimising the $1$--norm of $g(t)$.

\begin{thm}
	Suppose the maximum $\gamma > g(0)$. Then $\gamma$ is minimised if and only if:
	\[
		\Gamma := \int_{\R_+} g(t) \, dt
	\]%
	is minimised.
\end{thm}

\begin{pf}
	This follows by \Cref{lem:integrals}.
\end{pf}


\Cref{thm:fixed} gives conditions for $u$ to be optimal when either the input time $t'$ or duration $\tau$ is fixed. The input is optimal if the duration of the input is as short as possible so that the response is $\lambda$--incident. For example if the disturbance occurs before the input then the input would have a short duration. On the other hand should the input time occur before the disturbance then the duration, of the input, needs to be extended to prevent $g(t)$ falling below the minimum $\lambda$. Similarly when the duration is fixed the input time is constrained by the minimum value. As the magnitude of the input is fixed by the magnitude of the disturbance, the optimal duration and input time would be $\tau = 0$ and $t'= 0$ i.e. an impulse. As this would ensure that the response $g(t) \leq 0$ for all $t$. However such a duration and input time would, in general, result in the existence of a $t$ such that $g(t) < \lambda$.

\begin{thm}
	\labelx{thm:fixed}
	Fix $\lambda$ and suppose $U$ is adequate. Then:
	\begin{enumerate}
		\item{for fixed $t'$, $\gamma$ is minimised if and only if $\tau = \min\{\sigma : g(\sigma) \geq \lambda \wedge \exists t_{\min}; g(t_{\min}) = \lambda \}$.}
		\item{for fixed $\tau$, $\gamma$ is minimised if and only if $t' = \min\{s' : g(s') \geq \lambda \wedge \exists t_{\min}; g(t_{\min}) = \lambda \}$.}
	\end{enumerate}
\end{thm}

\begin{pf}
	Fix $w$ and $\lambda$.
	\begin{proofcase}
		Suppose $\tau > \sigma$ and let $u$ and $v$ be two adequate inputs with durations $\tau$ and $\sigma$ respectively but with the same input time $t'$. As $u$ and $v$ are adequate we have that:
		\[
			\int_{t'} ^s u <  \int_{t'} ^s v \leq U
		\]%
		for all $t' \leq s < t' + \tau$, where the strict inequality follows as the end point of the input $u$ is $'t  + \tau > t' + \sigma$ and $\V u \V_1  = \V v \V_1$. This holds only if $v(t) > u(t)$ for all $t \in [t', t'+\sigma]$. Thus by \Cref{lem:integrals} we have that $g(u) > g(v)$ for all $t > t'$.
	\end{proofcase}

	\begin{proofcase}
		This follows similarly by \Cref{lem:integrals}.
	\end{proofcase}
\end{pf}

\Cref{thm:main} provides the optimality conditions for the Magdelaine model for inputs $u$ of the form \eqref{eq:upulse}. Similarly to the results of \cite{town17} there are two conditions for optimality. In the first condition, should all minima occur prior to the global maximum of the response $g(t)$ then the optimal input is an input for which the duration $\tau = 0$. This case is similar to the optimality condition for the Bergman minimal model, derived by \cite{town17IFAC}, that the global maximum occurs between two global minima. However due to the requirement that $g(\infty) = 0$ and the instability of the equilibrium $g = 0$ in the Magdelaine model, there may not exist a second minimum of $g(t)$ which occurs after the maximum.
The second condition for optimality of an input to the Magdelaine model is identical to the condition for the Bergman minimal model found in \cite{town17IFAC} i.e. that the global minimum occurs between two global maxima. This is as the input is adequate and therefore guaranteed to return $g(t)$ to $0$.

\begin{thm}
	\labelx{thm:main}
	Fix $\lambda$ and suppose $U$ is adequate. 
	\begin{enumerate}
		\item{Suppose, for all $\tau$ and $t'$ that $\max \{t : g(t) = \lambda \} \leq \max \{ t : g(t) = \gamma\}$. Then $\gamma$ is minimised if and only if $\tau = 0$ and $t' = \min\{s' : g(s') \geq \lambda \wedge \exists t_{\min}; g(t_{\min}) = \lambda \}$.}
		\item{Suppose there exist $\tau$ and $t'$ such that $\max \{t : g(t) = \lambda \} \geq \max \{ t : g(t) = \gamma\}$. Then $\gamma$ is minimised if and only if there is $t_{\min} \in \argmin\{g(t)\}$ such that $\max_{t < t_{\min}} \{ g(t) \} =  \max_{t > t_{\min}} \{ g(t) \} = \gamma $.}
	\end{enumerate}
\end{thm}

\begin{pf}
	Fix $w$ and $\lambda$. Throughout this proof we say $g(u) > g(v)$ \emph{initially} if there exists $\varepsilon > 0$ such that $g(u) > g(v)$ for all $t \in (\min\{ s', t'\}, \min\{ s', t'\} + \varepsilon)$. 

	\begin{proofpart}
		Suppose $u$ is an input with duration $\tau = 0$ and $t'$ is such that $u$ is incident. Additionally, suppose there exists a distinct input $v \neq u$ such that $\gamma(v) < \gamma(u)$. In particular this implies that $g(v, t_{\max}) < g(u, t_{\max}) = \gamma(u)$. Additionally, as $\lambda$ is a fixed lower bound $g(v, t_{\min}) \geq g(u, t_{\min}) = \lambda$. 

		As $t_{\min} < t_{\max}$ this implies $t_g \in [t_{\min}, t_{\max})$, where $t_g$ is the intersection point of the responses $g(u)$ and $g(v)$. By \Cref{lem:intpoints,lem:integrals} this $t_g$ must be unique. Thus $g(v) > g(u)$ for all $t \in (t', t_g)$. This is true if and only if $u > v$ initially. This occurs if either $\sigma > \tau$ -- as the 1-norm of $u$ and $v$ are bounded -- or $s' > t'$. In all cases this implies:
		\[
			\int x(u) > \int x(v)
		\]%
		for all $t > t'$. Thus $g(u) < g(v)$ for all $t > t'$. Thus $g$ is minimised by $u$.

		Suppose $u$ is minimal but either $\tau > 0$ or $t'$ is such that $g(u) > \lambda$ for all $t$. In the latter case, by \Cref{thm:fixed} there exists input $v$ with the same duration $\tau$ as $u$ such that $\gamma(v) < \gamma(u)$. Thus, we may assume $t'$ is such that $u$ is incident. Suppose $v$ is an incident input with duration $\sigma < \tau$ and input time $s' > t'$. This implies $g(v) > g(u)$ initially and therefore by \Cref{lem:integrals} that $t_g \in [t_{\min}, t_{\max}]$. Thus $\gamma(v) < \gamma(u)$.  
	\end{proofpart}

	\begin{proofpart}
		Suppose $u$ is an incident input, with duration $\tau$ and input time $t'$, for which $\max \{t : g(t) = \lambda \} \geq \max \{ t : g(t) = \gamma\}$ and $v$ is an incident input, with duration $\sigma$ and input time $s'$, whose response is as in Part 1 of this Theorem. We say $v$ is a type 1 input and $u$ is a type 2 input. Define the sequence of adequate inputs $(o)_{i = 0} ^ \infty$ with durations $\alpha_i$ where:
		\[
			\alpha_{-1} := \tau, \quad \alpha_{-2} := \sigma, \quad \alpha_0 := \frac{\tau + \sigma}{2}
		\]%
		and:
		\[
			\alpha_i := \frac{\alpha_{i-1} (i-1) + \omega}{i}
		\]%
		where:
		\[
			\omega := 
			\begin{cases}
				\tau,   &  \alpha_{i-1} \text{ is type 1} \\
				\sigma, &  \alpha_{i-1} \text{ is type 2}
			\end{cases}
		\]%
		
		We partition $(o_i)$ into the two subsequences: $(p_j)$ and $(q_m)$ where $o_i \in (p_j)$ if $o_i$ is type 1 and $o_i \in q_m$ if $o_i$ is type 2. We now show that $\gamma(p_i) < \gamma(p_{i-1})$ and $\gamma(q_i) < \gamma(q_{i-1})$. Take $p_{i-1}$ and $p_i$ two subsequent elements of the sequence $(p_k)$ which have durations $\beta_{i-1}$ and $\beta_i$, respectively. By construction of the sequence $\beta_i < \beta_{i-1}$. Thus $p_{i}$ is nested in $p_{i-1}$ and therefore $g(p_{i-1}) < g(p_{i})$ initially. As $p_i$ and $p_{i-1}$ are type 1 and nested we have that $t_g \in [t_{\min, i-1}, t_{\max, i-1})$ as if it where not either there would exist $t$ such that $g(p_i,t) < \lambda$ or $p_i$ would not be type 1 as its maximum would occur before its minimum. We argue similarly to show  $\gamma(q_i) < \gamma(q_{i-1})$. 

		Let $L$ be the set of indexes for which $o_l$ and $o_{l-1}$ are inputs of different types. Without loss of generality, we assume $o_l$ is type 2 and $o_{l-1}$ is type 1. The difference in durations:
		\begin{align*}
			\alpha_{l} - \alpha_{l-1} & = \frac{\tau - \alpha_{l-1}}{l} \\
														 & = \l ( \frac{1}{l} \r ) \l ( \frac{1}{l-1} \r )  (\tau(l-1) + \tau (k + 1 - l) - \sigma k) \\
								       & = \l ( \frac{1}{l} \r ) \l ( \frac{k}{l-1} \r ) ( \tau - \sigma) \\
								& \leq  \l ( \frac{1}{l} \r ) (\tau - \sigma)
		\end{align*}%
		where $k \leq l - 1$ is some natural number. A similar expression holds should the types of $o_l$ and $o_{l-1}$ be reversed. Thus $\alpha_l - \alpha_{l-1} \to 0$ as $l \to \infty$. 

		As $g$ and therefore $\gamma$ are continuous functions of the duration we have that for all $\varepsilon > 0$ there exists $N$ such that $|\gamma(p_j) - \gamma(q_k)| < \varepsilon$ for any $j,k \geq N$. This occurs only if:
		\[
			\lim_{j \to \infty} \gamma(p_j) = \lim_{k \to \infty} \gamma(q_k)
		\]%
		Indeed for any $t$ we have:
		\[
			\lim_{j \to \infty} g(p_j) = \lim_{k \to \infty} g(q_k)
		\]%
		We need only consider the shapes of the sequence of responses $(p_i)$ to determine the shape of the response to the limit $o$. By the above we have that:
		\[
			g(p_{i-1}) < g(p_i) | _{B_i}
		\]%
		where $B_i := [t'_i , t_{g,i}]^c$. By the assumption that $p_i$ is type 1 we have that $ m_i := \max_{t < t_{\min}} \{ g(p_i)\} < \gamma(p_i)$. Thus $(m_i)$ is a monotone increasing sequence, bounded above by $\gamma(p_i)$ for each $i$ and has $\sup \{m_i\} = \gamma(p)$, where $p := \lim p_i$. Thus $m := \lim m_i  = \gamma(p)$ i.e. the global maximum before the global minimum equals the global maximum after the global minimum.

		Now let $g$ be a response to an input $u$ as per the statement of the Theorem. Suppose there exists an input $v \neq u$ such that $\gamma(v) < \gamma(u)$. If $g(v) > g(u)$ initially there must exist $t_g< t_{\max, 1}$ -- the point at which first maximum of $g(u)$ occurs. As there is at most one $t_g$ and $t_{\min} > t_{\max,1}$ this implies there exists $s$ such that $g(v,s) < \lambda$. Instead suppose $g(v) < g(u)$ initially. By the lower bound constraint there must exist $t_g \geq t_{\min}$. Again as $t_g$ is unique and there exists $t_{\max} > t_{\min}$ we have that there must exist $s$ such that $g(v,s) > \gamma(u)$. Hence no such $v$ exists.
	\end{proofpart}
\end{pf}

\begin{cor}
	Suppose the conditions of \Cref{thm:main} part 2 are met. Then there exists an input $u(t', \tau)$ which produces the minimised response.
\end{cor}

\section{Numerical Example}

\begin{figure}[h]
\begin{center}
\includegraphics[scale=0.63]{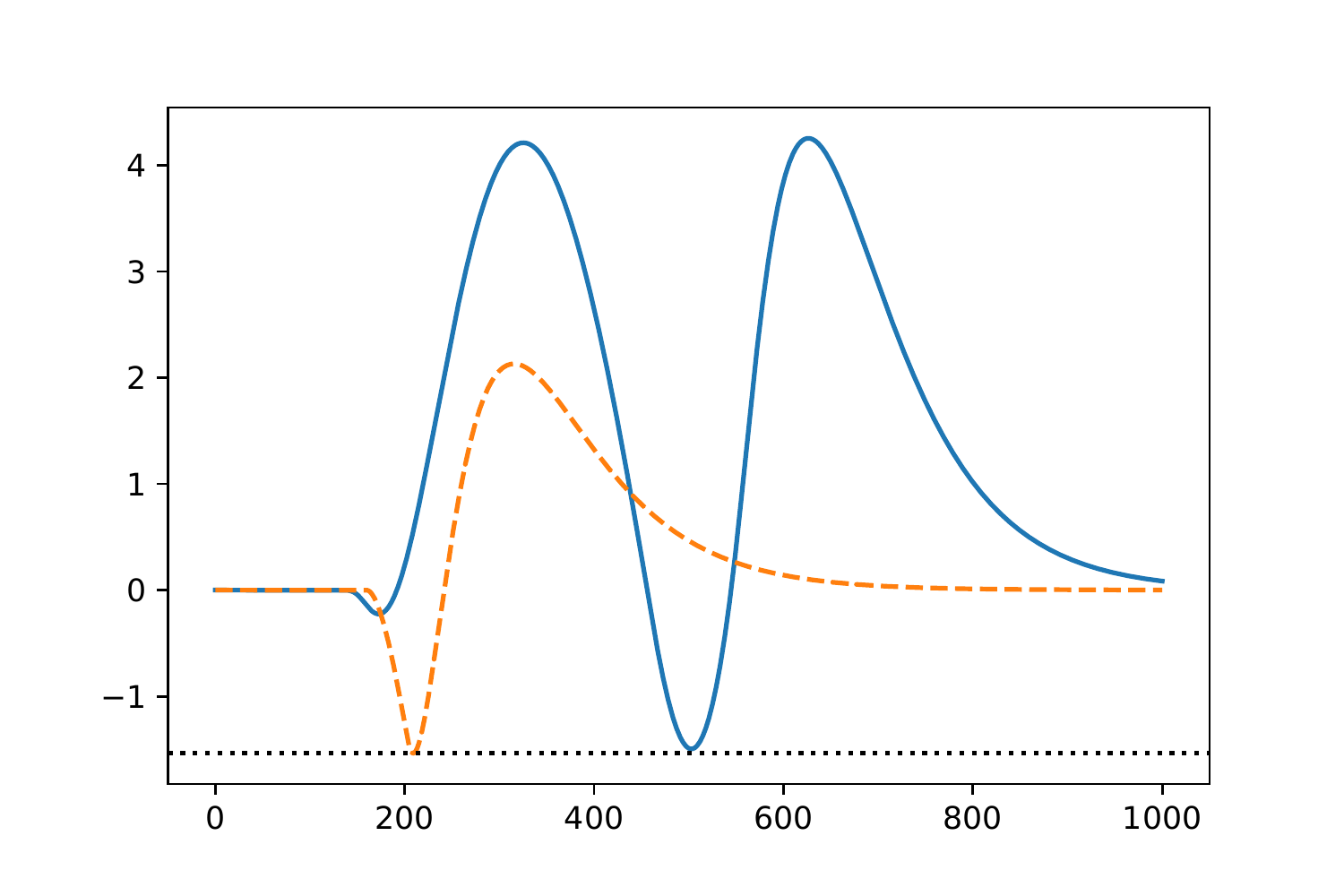}
\end{center}
\caption{The response of the Magdelaine model to the optimal inputs for a long duration disturbance (blue response) and a short disturbance (dashed orange response).}
\label{fig:nicgs}
\end{figure}

\Cref{fig:nicgs} shows the plasma glucose concentration of the Magdelaine model in response to the optimal pulse input for two different disturbances. These responses are normalised so the steady-state concentration is at $0 \, \mathrm{mmolL^{-1}}$.

For both responses the fixed lower bound was set as $\lambda =  -1.5$ which corresponds to a lower bound of $4.0 \, \mathrm{mmolL^{-1}}$ in the non-normalised model. The blue response is an example of the second optimality condition given in \Cref{thm:main}. As two equal maxima occur about the global minimum.  The disturbance is:
\[
d(t) := \chi_{[150, 400]}
\]%
The optimal bolus for this disturbance is $0.075 \chi_{[133, 448]}$ i.e. a pulse with input time $t' = 133$, duration $\tau = 315$ and magnitude $\hat u = 0.075$.

The dashed orange response is an example of the first optimality condition given in \Cref{thm:main}. As the global minimum occurs before the global maximum and the duration of the input is $\tau = 0$. It is a response of the disturbance:
\[
d(t) := 20 \chi_{[200, 202]}
\]%
The optimal input is:
\[
u(t) := 2.36 \delta(t - 158)
\]%
i.e. a pulse with duration $\tau = 0$ and input time $t' = 158$.

\section{Bounded Inputs to the Bergman Minimal Model}

The Bergman Minimal model (\cite{berg05}, \cite{kand09}), is a \emph{non-linear} continuous-time model of glucose and insulin dynamics in type one diabetes which is used as the basis of more complicated models such as those of \cite{fabi06} and \cite{kand09}. In contrast with the Magdelaine model, the Bergman model depends recursively on the current glucose state. The model is comprised of a set of first order linear ordinary differential equations which govern the subcutaneous, plasma and interstitial concentrations and effectiveness of insulin, denoted by $z, y$ and $x$ respectively:
\begin{align}
\labelx{eq:optbergmod}
\begin{pmatrix}
\dot x\\
\dot y\\
\dot z
\end{pmatrix}
=
\begin{pmatrix}
-a & ab & 0 \\
0  & -c & c\\
0  & 0  & -d
\end{pmatrix}
\begin{pmatrix}
x\\
y\\
z
\end{pmatrix}
+
\begin{pmatrix}
0\\
0\\
dk
\end{pmatrix}
u
\end{align}%
and a non-linear ordinary differential equation which governs the plasma glucose concentration $g(t)$:
\begin{align}
\labelx{eq:optbergmnonlin}
\dot g(t) = -(x(t) + G)g(t) + w(t)
\end{align}

In \eqref{eq:optbergmod} the parameters $a, b, c, d$ and $k$ are positive time constants which control the rate of transfer of insulin between the states $z, y$ and $x$. The constant $G > 0$ in \eqref{eq:optbergmnonlin} represents  insulin independent glucose uptake or loss e.g. via renal excretion.

%
As mentioned above, \cite{town17} and \cite{town17IFAC} characterised the optimality of pulse inputs to the Bergman minimal model, in terms of the glucose response, for any given bounded disturbance $w$. As any positive plasma glucose concentration less than an upper bound, determined by the constant $G$, is an asymptotically stable equilibrium determined by the basal input $\overline u$, the plasma glucose concentration will always return to steady-state independent of the total amount of bolus insulin delivered. Thus the $1$-norm of inputs to the Bergman minimal model is not constrained by the requirement to return to steady-state as it is in the Magdelaine model. 

In \cite{town17IFAC} the duration $\tau$ of the input $u$ is fixed and the optimal input time $T$ and magnitude $\hat u$ is found. This is extended in \cite{town17} to optimise the input duration. Thus the pulse input which minimises the magnitude of $g$ whilst remaining above a lower bound $\lambda$ is characterised in terms of the response of $g$ with respect to $u$. 

\begin{figure}[h]
\begin{center}
\includegraphics[scale = 0.6]{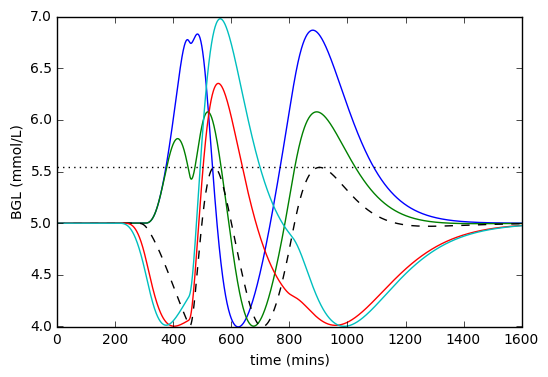}
\caption{Optimal pulse inputs to the Bergman minimal model for a variety of fixed durations \citep{town17}.}
\label{fig:exampleoptberg}
\end{center}
\end{figure}

We will say an input is \emph{optimal in the sense of} \cite{town17} if the response to the input satisfies the optimality conditions derived in \cite{town17} i.e. if the maxima and minima of the response are iterlaced. An example of inputs which are optimal in the sense of \cite{town17} and \cite{town17IFAC} are shown in \Cref{fig:exampleoptberg}. 

The dark blue and green responses in \Cref{fig:exampleoptberg} are optimal in the sense of \cite{town17IFAC} as the maxima on either side of the minimum are equal whereas the light blue and red responses are optimal as the global maximum occurs between two global minima.

\cite{town17} proved that it is possible to further optimise the response by optimising the duration of the bolus input. The optimal response is given by the dashed black response in \Cref{fig:exampleoptberg}. 

Here we additionally constrain the $1$-norm of the bolus input to the Bergman minimal model. This constraint on the total amount of bolus insulin delivered may be used as a more feasible constraint to avoid the potential risk of over bolusing insulin resulting in hypoglycaemia than specifying a lower bound above the glycaemic threshold and is more robust to errors in estimation of the disturbance $w$.

However, this constraint alters the optimality conditions of \cite{town17}. Given a specified lower bound, $\lambda$, there could exist a response which does not attain the specified minimum $\lambda$ yet has a lower maximum than a response which does obtain the specified minimum.

For a given disturbance $w$ and fixed lower bound $\lambda$ we will take the \emph{required} bolus amount to be:
\[
U :=	\l \V \hat u \chi_A \r \V := \int_{[0, \infty)} \hat u \chi_A \, dt
\]%
to be the amount so that the response is optimal in the sense of \cite{town17} and consider the optimality of inputs which are less than this amount. In \Cref{thm:bergext} we suppose $u$ is a pulse input of the form \eqref{eq:u} to the Bergman minimal model for which the bolus is less than the required amount. 

Throughout the remainder of this section we fix $\lambda$ and let $w$ be a bounded positive disturbance with a required bolus amount $U$. We also take $u$ and $v$ to be pulse inputs of the form \eqref{eq:upulse} -- with input times $t'$ and $s'$ and durations $\tau$ and $\sigma$, respectively. Furthermore we set the bolus amounts of $u$ and $v$ to be identical i.e. $\V \hat u \chi_{[t', t' + \tau]} \V_1  = \V \hat v \chi_{[s', s' + \sigma]} \V_1 < U$ and take $\overline u = \overline v$. As the global minimum attained by $g$ in response to the input $u$ is no longer fixed to be $\lambda$, we define $\lambda(u) := \min \{ g(u) \}$.

\begin{thm}
\labelx{thm:bergext}
Suppose for all minima $t_{\min}$ of the response $g$ that $\max_{t < t_{\min}} \{ g(t) \} \ne \max_{t > t_{\min}} \{ g(t) \}$ for all pulse inputs $u$ such that $\V \hat u \chi_A \V_1 \leq U$. Then $\gamma(u) < \gamma(v)$ if and only if $\tau < \sigma$.
\end{thm}

\begin{pf}
Let $v$ be an input with duration $\sigma > 0$ and $u$ an input with duration $\tau < \sigma$. 

Should $t' = s'$ then we have that $g(u) < g(v)$ initially. Similarly if $t' \geq s_{\max}$ then we know $g(u) > g(v)$ for all $t \leq s_{\max}$ which implies $\gamma(u) > \gamma(v)$. As $g(u)$ is a continuous function of $t'$ and there are $t'$ as above, there must exist a $t' \in (s', s_{\max})$ and $t_g \in (t', s_{\max})$ such that $g(u) > g(v)$ for all $t < t_g$, $g(u) = g(v)$ when $t = t_g$ and $g(v) > g(u)$ for all $t > t_g$. Thus $g(u) < g(v)$ for all $t \geq s_{\max}$ which implies $\gamma(u) < \gamma(v)$.
\end{pf}

\begin{cor}
\labelx{cor:corextend}
Suppose the maxima of the response to the required bolus occur between two global minima. Then any bolus less than the required bolus is optimal if and only if $\tau = 0$.
\end{cor}

\begin{pf}
According to the results of \cite{town17}, the duration of the required bolus is $\tau = 0$. As the model is monotonic in the input $u$ if $\V u \V_1 < U$ then the reponse $g(u) > \lambda$ for all $t$. Additionally \Cref{thm:bergext} implies that the maximum of the response $\max \{ g(u) \}$ is minimised when the duration $\tau$ is minimised. Thus the duration of the input $u$ must be $0$. 
\end{pf}

\section{Example of Constrained Optimality Condition}

\begin{figure}[H]
\begin{center}
\includegraphics[scale = 0.65]{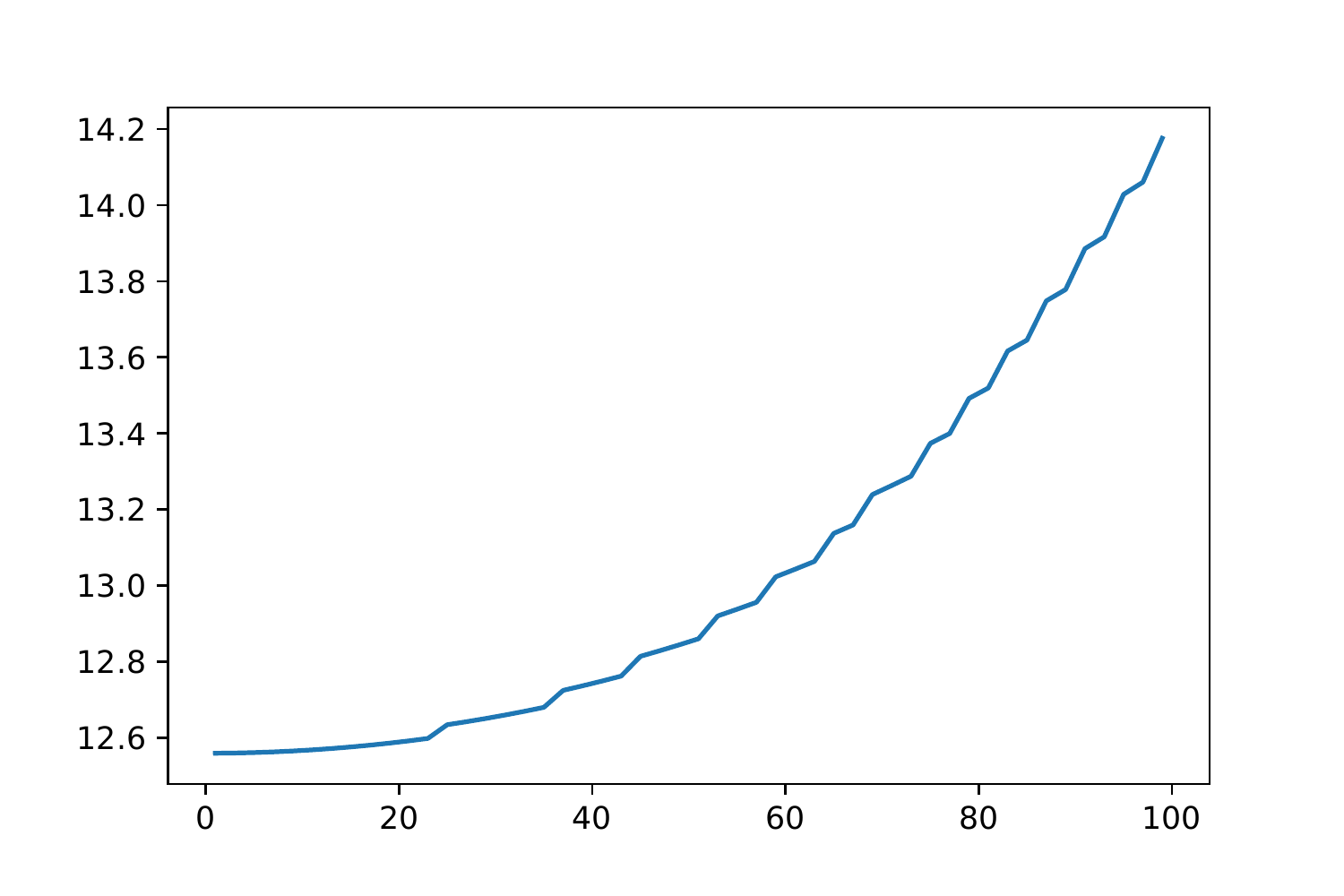}
\end{center}
\caption{The maximum of the response of the Bergman model to an input as a function of the duration of the input where the amount of the input is fixed to be a number less than the required bolus for the disturbance.}
\label{fig:bergdurs}
\end{figure}

The example presented in \Cref{fig:bergdurs} shows the maximum of the response of the Bergman minimal model to constrained inputs of various durations -- where the disturbance $w(t) := 263^{-1}f_1(t) + 1.0$ where $f_1$ is the solution to:
	
\begin{align}
\label{eq:genw}
\begin{pmatrix}
\dot f_1 \\
\dot f_2
\end{pmatrix} & = 
\frac{1}{60}
\begin{pmatrix}
-1 & 1\\
0 & -1
\end{pmatrix}
\begin{pmatrix}
 f_1 \\
 f_2
\end{pmatrix} 
 +
 \begin{pmatrix}
 0 \\
 4
 \end{pmatrix}
 \chi_{[200, 202]}(t)
\end{align}

The unconstrained optimal pulse input, i.e. the input which is optimal in the sense of \cite{town17}, is:
\[
u(t) :=  35.15 \delta (t - 175)
\]%
The response to this input has a global maximum of $8.5 \, \mathrm{mmolL^{-1}}$ which occurs beteen two global minima. 
The inputs in the example presented in \Cref{fig:bergdurs} are of the form:
\[
u(t) :=  \l ( \frac{20}{\tau} \r ) \chi_{[t', t' + \tau]}
\]%
This constrains the total amount of bolus insulin to be $20 < U = 35.15$ -- which is the required bolus. The input time $t'$, for each duration $\tau$, is taken to be:
\[
t' := \argmin \l \{ g(u(t', \tau)) : \lambda (u) \geq 4.0 \r \}
\]%
i.e. the input time which minimises the maximum plasma glucose concentration. As shown the lowest maximum plasma glucose concentration occurs when $\tau = 0$. The jaggedness of the plot is an artefact of the numerical precision of the simulation in which the input time $t'$ was restricted to be an integer.

\section{Conclusions}

We have characterised the optimality of bolus inputs to the Magdelaine and Bergman models of type one diabetes when the total volume of insulin is constrained. This constraint arises from the structure of the Magdelaine model as it is necessary for inputs to meet this constraint to return the plasma glucose concentration to steady-state. We have proven that an input is optimal when the minimum of the plasma glucose response occurs prior to the maximum or if the minimum occurs between two equal maxima. Any further attempt to lower peak plasma glucose concentration will result in the plasma glucose concentration dropping below the fixed lower bound i.e. hypoglycaemia.

For the Bergman model the input which minimises the maximum plasma glucose concentration does not necessarily attain the lower bound. This differs from the results of \cite{town17IFAC} and \cite{town17} in which the volume of insulin delivered was not constrained. This suggests that the duration and timing of a bolus input are as significant as the total volume delivered.

Further work will focus on characterising the optimality of constrained inputs to the Bergman minimal model when there is an input such that the maxima on either side of the global minimum are equal. This case is not covered by \Cref{thm:bergext}.

It is also of interest to investigate optimality conditions for both the Bergman and Magdelaine models when it is possible to lower the basal insulin flow. For example to set $\overline u = 0$ on some bounded interval. In the Magdelaine model we expect setting $\overline u = 0$ on some interval will allow the results of \cite{town17} to apply directly.


\bibliography{references}

\end{document}